\documentclass[mathpazo]{aamm}
\renewcommand\thisnumber{x}
\renewcommand\thisyear {200x}
\renewcommand\thismonth{xxx}
\renewcommand\thisvolume{xx}
\usepackage{color}
\usepackage{algorithm}
\usepackage[noend]{algpseudocode}

\definecolor{manuelc}{rgb}{1, 0, 0}

\newcommand{\calL}{{\mathcal L}}

\newcommand{\N}{{\cal N}}

\newcommand{\D}{{\cal D}}

\newcommand\rev[1]{\textcolor{red}{#1}}

\newcommand{\bs}[1]{\boldsymbol{#1}}
\newcommand{\bmu}{\bs{\mu}}

\begin{document}


\markboth{Ji, Chen, Xu}{RBM for Nonlinear PB}
\title{A reduced basis method for the nonlinear Poisson-Boltzmann equation}



%
%
%

\author[Lijie Ji, Yanlai Chen and Zhenli Xu]{Lijie Ji\affil{1},   Yanlai Chen\affil{2}\comma\corrauth  and Zhenli Xu\affil{3}}
\address{\affilnum{1}\ School of Mathematical Sciences, Shanghai Jiao Tong University, Shanghai 200240, China. Email: {\tt sjtujidreamer@sjtu.edu.cn}.\\
\affilnum{2}\ Department of Mathematics, University of Massachusetts Dartmouth, 285 Old Westport Road, North
Dartmouth, MA 02747, USA. Email: {\tt yanlai.chen@umassd.edu}.\\
\affilnum{3}\ School of Mathematical Sciences, Institute of Natural Sciences,
and Key Laboratory of Scientific and Engineering Computing (Ministry of Education),
  Shanghai Jiao Tong University, Shanghai 200240, China. Email: {\tt xuzl@sjtu.edu.cn}.
}
\begin{abstract}
In numerical simulations of many charged systems at the micro/nano scale, a common theme is the repeated solution of the Poisson-Boltzmann equation. This task proves challenging, if not entirely infeasible, largely due to the nonlinearity of the equation and the high dimensionality of the physical and parametric domains with the latter emulating the system configuration. In this paper, we for the first time adapt a mathematically rigorous and computationally efficient model order reduction paradigm, the so-called reduced basis method (RBM), to mitigate this challenge. We adopt a finite difference method as the mandatory underlying scheme to produce the {\em truth approximations} of the RBM upon which the fast algorithm is built and its performance is measured against. Numerical tests presented in this paper demonstrate the high efficiency and accuracy of the fast algorithm, the reliability of its error estimation, as well as its capability in effectively capturing the boundary layer.

\end{abstract}



\maketitle
\graphicspath{{Figs/}}

\section{Introduction}

{Fast numerical algorithms for solving} parameterized {partial differential equations (PDEs) have attracted wide-spread
interest in recent years}, particularly in engineering {applications due to} many control, optimization and design
problems {requiring repeated simulation of} certain parameterized PDEs. Traditional numerical methods solve the equation for each necessary parameter value and thus {obtaining the solution ensemble}
for the whole parameter space is potentially time-consuming {if not entirely infeasible.} This is an especially onerous task if the physical and/or parametric domain are of high-dimension. {It is therefore imperative to design} efficient and accurate {reduced order modeling techniques} for these scenarios capable of realizing negligible marginal (i.e. per parameter value) computational cost. 
The reduced basis method (RBM) provides a rigorous and highly efficient platform to achieve this exact goal. It was first introduced for nonlinear structure problem \cite{Almroth1978, noor1979reduced} in 1970s  and has been later
analyzed and extended to solve many problems such as linear evolutionary equation \cite{HaasdonkOhlberger}, viscous Burgers equation \cite{veroy2003reduced}, Navier-Stokes equations\cite{deparis2009reduced}, and harmonic Maxwell's equation \cite{chen2010certified, chen2012certified} just to name a few. 
Interested readers are referred to the review papers \cite{Rozza2008, Haasdonk2017Review} and recent monographs \cite{Quarteroni2015, HesthavenRozzaStammBook} for a systematic description of the RBM.

One such parametric scenario we are concerned in this paper is the simulation of the electrostatic interaction which is essential for many systems in  physical, biological and materials sciences \cite{FPP+:RMP:2010,WK+:N:2011,Levin:RPP:2002} at the nano/micro scale. These include, for example, biopolymers, colloidal suspensions, and electrochemical energy devices. 
The Poisson-Boltzmann (PB) theory \cite{Gouy:JP:1910,Chapman:PM:1913,Baker:COSB:2005,FBM:JMB:2002} plays a fundamental role in understanding the
electrostatic phenomenon in such systems. It subjects the electric potential of a charged system at the equilibrium state
to a nonlinear elliptic equation with the the Boltzmann distribution for the ionic densities.
The numerical solution of the PB equation has been widely studied in literature \cite{Baker:ME:2004a,LZH+:CCP:2008}, and the numerical
solvers are implemented in many popular software packages such as Delphi and APBS for practical simulations.
However, one often needs to solve the PB equation repeatedly to determine certain physical quantities of interest (QoI) which are usually dependent on a wide range of {parameters} delineating e.g. the boundary voltage, the
geometric length, and the Debye length. {Particular examples of such QoIs} include the
electrochemical capacitance, the current-voltage relation, and the free-energy calculation etc.

In this work, we propose a reduced basis method for the {parameterized} nonlinear PB equation.
Model order reduction for nonlinear equations is often realized by linearization techniques \cite{Weile:IEEE:2001} or polynomial approximations \cite{Phillips2003}, among others. One frequently-used tool  is the empirical interpolation method (EIM) \cite{Barrault2004, grepl2007efficient} which is crucial to facilitate the offline-online decomposition, a hallmark feature of RBM to realize the negligible marginal computational cost. This paper extends the RBM
for the nonlinear PB equation by approximating the nonlinear exponential term with a Taylor expansion form \cite{Shestakov2002}. This  leads to a linear equation in each calculation step.
Realizing a partial offline-online decomposition, the method promises high accuracy due to the avoidance of the EIM error.
It is noted that this work only focuses on the mean-field PB equation which is limited to describe phenomena when many-body interactions are important. The extension of our work to the modified PB equations such as those including correlation and steric effects {(see, e.g., \cite{LJX:SIAP:2018} and references therein)} is 
of great interest due to the complex electrostatic phenomena they model and the drastically different nonlinearity contained therein. The successful application of the RBM will be reported in the future. We also note that Kweyu {\it et al.} \cite{Kweyu2016, Kweyu2017} has recently extended the RBM to the linearized PB equation for electrostatic solvation calculations of the biomolecules. To the best of our knowledge, this paper is the first attempt of solving the fully nonlinear equation (with rapid nonlinearity) by the RBM.

The paper is organized as follows. In section \ref{sec:bg}, we introduce the basic RB algorithm.
Detailed description of the  PB model, the FDM scheme used to obtain our {\em truth approximation}, how we apply RB to the nonlinear PB equation and its computational analysis are provided in section \ref{sec:alg}.
In section \ref{sec:result}, we show numerical results in both one and two physical dimensional spaces to demonstrate the accuracy and efficiency of our reduced model.
Finally, concluding remarks are drawn in section \ref{sec:conclusion}.

\section{Overview of the reduced basis method}
\label{sec:bg}

The reduced basis method is a fast algorithm for computing a certified surrogate to the highly accurate but potentially expensive numerical solution (termed {\em truth approximation} in the RB context)  of a system dependent on a $P$-dimensional parameter
\[
\bmu \in \D \subset {\mathbb R}^P.
\]
It is particularly useful for the many-query or real-time simulation context where an initial investment may pay off through repeated simulations with significantly less (at times negligible) marginal cost at a later stage. 
An essential tool is the offline-online decomposition process. The offline phase is devoted to construct the RB space, denoted by $W_N$ (with $N$ being its dimension, and usually much smaller than the number of degrees of freedom for the truth approximation.) 
During the online stage, a RB approximation for any given parameter value $\bmu$ in the prescribed domain $\D$ is sought from the space $W_N$.

We note that the online solver (of various dimensions) is invoked repeatedly offline to construct the RB space $W_N$ through a greedy algorithm. We therefore present here the crucial online solver for a linear PDE, and postpone the construction of the RB space until when we describe the RBM for the nonlinear PB equation for completeness of that section. 
Indeed, consider a linear elliptic PDE, $L(\bmu)u(x,y)=f(x,y)$, with the operator (and/or the right hand side) parameterized by $\bmu$.
Let $\N$ be the number of degrees of freedom for a well-defined and accurate numerical scheme (termed {\em truth solver} in the RB context) discretizing this equation.  The numerical approximation $u^\N(\bmu)$ is the solution of
\begin{equation}
\label{eq:pdesystem}
L_\N(\bmu) u^\N = f,
\end{equation}
which can be understood as deriving from a collocation formulation, $L_\N(\bmu) u^\N(x, y) = f(x, y)$ for certain $(x,y)$'s, or a Galerkin scheme, $a^\N(u^\N, v) = f(v) \mbox{ for } \forall \,\, v$.
A critical assumption for the operator $L_\N$ is that it is affine with respect to (the functions of) the parameter. That is, it can be written as,
\begin{equation}L_\N =\sum_{q = 1}^Q B_q(\bmu) L_{\N}^q
\end{equation}
where $L^q_{\N}$ is a parameter-independent operator and coefficient function $B_q$ depends on parameter $\bmu$.
The RBM is built upon this discrete solution and its accuracy is also measured against it. For that reason, the solution of Eq. \eqref{eq:pdesystem} is considered ``exact'' and thus called the {\em truth approximation}. For simplicity of exposition, we shall drop the superscript $\N$ in the remainder of the paper as we will not make any reference to the exact solution of the PDE.

The online process of the RBM is as follows. Assuming that we have identified $N$
parameter values $\{\bmu^1, \dots, \bmu^N\}$ and the corresponding truth approximations $u_n \equiv u^\N(\bmu^n), 1 \le n \le N$. 
With a slight abuse of notation, we don't differentiate these function and their discrete vectors. 
These vectors constitute the basis space of the RBM, written in the form of a matrix, $W_{ N} = [u_{1},  \ldots, u_{N}] \in \mathbb{R}^{\N \times N}$.
One expresses the RB approximation as a linear combination of the basis vectors. That is, we have
\begin{equation}
\label{eq3}
\widehat{u}(\bmu) = W_{ N}  c_N(\bmu) ,
\end{equation}
where $c_N (\bmu)\in \mathbb{R}^N$ is the RB coefficient vector. These coefficients are sought by satisfying the ansatz of the PDE. Therefore, we substitute this combination into equation \eqref{eq:pdesystem} to obtain a linear algebra system,
\begin{equation}
A(\bmu)c_N(\bmu) = f,\end{equation}
where $A(\bmu) = L_\N W_N$ is a $\N \times N$ matrix, $f \in \mathbb{R}^\N$ is a column vector,
and $\N \gg N$.

Coefficients $c_N (\bmu)$  can be obtained by the least-squares method leading to a Petrov-Galerkin approach.
One can also resort to a Galerkin approach (i.e. an orthogonal projection-based  order reduction technique \cite{Phillips2003})
to identify $c_N(\bmu)$ as a solution of the following ({\em reduced}) linear system:
\begin{equation}
\label{eq2}
B(\bmu) c_N(\bmu) = f_N.
\end{equation}
Here $B(\bmu)$ is the RB matrix of dimension $N \times N$ and $f_N$ is the RB vector of dimension $N$, which are expressed as follows,
\begin{equation}
\label{eq2b}
B(\bmu)  = W_{N}^T L_\N (\bmu) W_{ N}, \qquad  f_N(\bmu) = W_{ N}^T f .
\end{equation}
Obviously, this is an energy projection into the RB space ${\rm span}\{u_n:1 \le n \le N\}$. Solving equation \eqref{eq2} is much cheaper than solving equation \eqref{eq:pdesystem} and system \eqref{eq2}
is an order reduction in comparison to system \eqref{eq:pdesystem}.

\section{Poisson-Boltzmann equation and its reduced model}
\label{sec:alg}
\subsection{Poisson-Boltzmann model and the truth approximation}
\label{sec:bg_PB}

Poisson-Boltzmann  is a mean-field theory describing the equilibrium distribution of charged systems \cite{  Davis1990,Levitt:Biophy:1985,Jordan:1998,Weetman:1997},
which has been widely used in biomolecular solvation, microfluidic devices, and charged soft materials.
Typically, one considers a symmetric binary electrolyte between two parallel electrodes positioned 
at locations marked by $X=\pm L_X$ with the extremes of the other direction marked $Y = \pm L_Y$.
The PB equation for the electric potential $\phi$ is
written as,
\begin{align}
\label{eq:pb}
\nabla \cdot \left ( \varepsilon \nabla \phi \right )= 2zec_0\sinh(\beta ze\phi)-\rho_f,
\end{align}
where $\varepsilon$ is the dielectric permittivity, $\pm ze$ is the charge of an cation or anion,
$\beta$ is the inverse thermal energy, $c_0$ is the bulk concentration, and $\rho_f$ is the density of the fixed charge.

Without loss of generality, we let the computational domain be a 2D square by setting $L_X=L_Y=L$. We further set $x=X/L$, $y=Y/L$, $\varepsilon=\varepsilon_W$
being a constant, $\Phi=z e \beta\phi$, and $g=\rho_f/(2zec_0)$. Then Eq. \eqref{eq:pb} becomes the following dimensionless PB equation,
\begin{subequations}
\label{eq:PB}
\begin{align}
D\nabla^2 \Phi  = \sinh \Phi + g(x,y),\\
\intertext{where $D = {(\ell_D/L)}^2$, $\ell_D=1/\sqrt{8\pi\ell_Bz^2 c_0}$ is the Debye screening length, and
$\ell_B =\beta e^2 /(4\pi \varepsilon_W)$ is the Bjerrum length in water solvent. The Laplace operator
$\nabla^2$ is with respect to the new coordinates $(x,y)$, and the computational domain in $(x,y)$ is now  $\Omega=[-1,1]^2$. We introduce the following boundary conditions,}
\label{eq:PB_bc1}
\Phi(x=\pm 1, y)  = \pm V,  \\
\label{eq:PB_bc2}
\partial_y \Phi(x, y=\pm 1)  = 0.
\end{align}
\end{subequations}
Here, the first boundary condition represents the fixed boundary voltages on the left and right
electrodes, and the second one characterizes a state of low dielectric permittivity at top and bottom boundaries $y=\pm 1$. Equation \eqref{eq:PB} is thus the PB equation parametrized by
\[
\bmu := [D, V],
\]
a vector-valued parameter. We intend to devise a RBM for its rapid resolution for scenarios when it needs to be solved repeatedly for a wide range of $\bmu$ values.

\subsubsection{The truth solver}
Before the discussion of the RBM, let us describe a finite-difference solver for the nonlinear PB equation.
Numerical methods for nonlinear PB equations have been widely studied \cite{LZH+:CCP:2008,LuZhangMcCammon_JCP05,Boschitsch2004,BHW:JCC:2000,Holst2007}.
In this paper, we use the similar technique as Shestakov {et al.} \cite{Shestakov2002}, which transforms
the nonlinear equation into a linear one at each iteration by truncating the Taylor series.
Let $\Phi_m$ be the approximate solution at the $m$th iterative step, then for the solution at the
$(m+1)$th step, the nonlinear term $\sinh \Phi_{m+1}$ is approximated by,
\begin{equation}
\sinh \Phi_{m+1} \approx \sinh \Phi_{m} + (\cosh \Phi_m) (\Phi_{m+1}-\Phi_m)
\end{equation}
 and the PB equation then becomes,
\begin{equation}
- D \nabla^2 \Phi_{m+1} + (\cosh\Phi_m) \Phi_{m+1} = (\cosh\Phi_m) \Phi_m -\sinh \Phi_m.
\label{linearpb}
\end{equation}
We then use the second-order five-point central difference scheme to approximate $\nabla^2\Phi$, leading to a
linear system for $\Phi_{m+1}$.
To describe this system, we denote by $\mathcal{N}$ the total number of grid points discretizing the physical domain. The number of free nodes, i.e. those in the interior of the domain and on boundaries $y=\pm 1$  is denoted by $\mathcal{N}_0$. This means that there are $\mathcal{N}-\mathcal{N}_0$ nodes on boundaries $x=\pm 1$ for which the corresponding potential values are specified.
Let $\mathcal{L}_\mathcal{N}(\bmu; \Phi)$ be the discretized operator for approximating the left hand side of
equation \eqref{linearpb} and Neumann boundary condition \eqref{eq:PB_bc2}. Let $\vec{\Phi}$ be the $(\N_0 \times 1)$ vector representing the discretized function $\Phi(x, y)$. The numerical scheme can then be written as
\begin{equation}
\mathcal{L}_{\mathcal{N}}(\bmu; \Phi_m) \vec{\Phi}_{m+1} = \vec{F}(\Phi_m),
\label{Operator}
\end{equation}
for $m =0, 1, \cdots$. Here $\vec{F}$ discretizes the right hand side of equation \eqref{linearpb} and incorporates the Dirichlet boundary condition \eqref{eq:PB_bc1}.  The iterative algorithm for solving \eqref{Operator} is summarized in Algorithm~\eqref{alg1} and the resulting solution $\Phi (\bmu)$ is called the ``truth'' approximation/solution corresponding to parameter $\bmu$ in the RBM framework.

\begin{algorithm}[h]
\caption{Iterative solver for the nonlinear PB equation}
\label{alg1}
\begin{algorithmic}[1]
\State Initialize potential distribution $\Phi_{0}$ and the  tolerance $\delta_0 = 1$.
\While {$ \delta_m > 10^{-11}$}
\State Solving the linear system of equations \eqref{Operator}; 
\State Set $\delta_{m+1} = ||\vec{\Phi}_{m+1} - \vec{\Phi}_{ m}||_{\infty} $;
\State Set $ m = m+1$;
\EndWhile
\State $\Phi (\bmu) =\Phi_m$.
\end{algorithmic}
\end{algorithm}

\subsubsection{The quantity of interest}

For most parametric systems, there are frequently quantities of interest which are nothing but functions of the parameter(s) describing the system. These QoIs are often calculated as functionals of the field variable, i.e. solution of the PDE modeling the system. Therefore, the efficient resolution of these field variables immediately leads to that of the QoI.

The electrochemical systems \cite{Conway::1999} of interest in this paper are no exceptions. Indeed, we are concerned with
the {\em total differential capacitance} of the symmetric electrolyte. It is defined as $C = C_L/2$, where $C_L$ is
the {differential capacitance} of the left electrode defined by
\[
C_L = \frac{d \sigma(V)}{d V} \quad \mbox{ with } \quad \sigma(V) = D \frac{\partial \overline{\Phi}}{\partial x} (x = -1),
\]
where $\sigma(V)$ is the {\em surface charge density}  at the left electrode. $\overline{\Phi}(x)$ is the average electric potential which is simply $\Phi(x)$ if the physical domain is one-dimensional and, for 2D, is nothing but
\[
\overline{\Phi} (x) =\int_{-1}^1 \Phi(x = -1, y)  dy.
\]

For the one-dimensional case, one can derive an explicit expression for the differential capacitance $C_L$
by solving the PB equation. Indeed, integrating Eq.~\eqref{eq:PB} from $x$ to $0$ gives
\begin{equation}
\frac{d \Phi}{d x}= -\frac{2}{\sqrt{D}}\sinh\left(\frac{\Phi}{2}\right),
\end{equation}
then one has, by utilizing the boundary condition, that $\sigma =-2 \sqrt{D}\sinh(-V/2)$, which means,
\begin{equation}
C_L = - \frac{d \sigma}{dV} = \sqrt{D} \cosh\left( \frac{V}{2}\right).
\label{eq:capacitance}
\end{equation}

\subsection{Reduced basis method for the Poisson-Boltzmann equation}
\label{sec:RBMonline}

As shown in the overview, the online procedure of the RBM algorithm is to find the coefficients of
the surrogate solution in the reduced basis space. Indeed, the $N$ dimensional coefficient vector $c$ is sought by
asking the resulting surrogate solution to satisfy the PDE \eqref{Operator} weakly in the RB space,
\begin{align}
A(\bmu; \widehat{\Phi}_{m}) c_{m+1} (\bmu) = W_N^T \vec{F}(\widehat{\Phi}_m), ~~m=0, 1, \cdots,
\label{eq:RBnlsys}
\end{align}
at every iteration. Note that $A(\bmu; \widehat{\Phi}_{m}) = W_N^T \mathcal{L}_\N(\bmu; \widehat{\Phi}_{m}) W_N$ depends on the current iterate of the RB solution
\begin{equation}
\widehat{\Phi}_m(\bmu) = W_N c_m(\bmu).
\end{equation}
Thus, we have to rely on an online iterative solver as well.
This iterative procedure for solving the coefficient $c(\bmu)$ is summarized in Algorithm \ref{alg2} for any given parameter $\bmu$
%
with the final RB surrogate approximation denoted by $c(\bmu)$.

\begin{algorithm}[h]
\caption{RB approximation for nonlinear PB equation, $c(\bmu) = $RBM$\_$PB$(W_N, \bmu)$ }
\label{alg2}
\begin{algorithmic}[1]
\State Initialize the potential distribution  $\Phi_{0}$ and the  tolerance $\delta_0 = 1$.
\While {$ \delta_j > 10^{-8}$}
                 \State Form the coefficient matrix $A$ and $W_N^T\vec{F}$ at each $j$th iteration.
               \State Solve for $c_{j+1}$ from \eqref{eq:RBnlsys}.
               \State $\Phi_{j+1} = W_N c_{j+1}$.
                 \State $\delta_{j+1} = ||\Phi_{ j}  - \Phi_{j+1}||_{\infty} $
                 \State j = j+1
\EndWhile
\State $c(\bmu) = c_{j}$.
\end{algorithmic}
\end{algorithm}
The full offline algorithm for constructing RB basis space  $W_N$ is realized with  standard greedy algorithm \cite{Rozza2008, Quarteroni2015, HesthavenRozzaStammBook}, which exploits a rigorous (albeit costly) {\em a posteriori} error estimator. Discretizing the parameter domain $\D$ by a sufficiently fine training set $\Xi_{\rm train}$, the greedy algorithm starts by selecting  the first  parameter $\bmu^1$ randomly from $\Xi_{\rm train}$ and obtaining its corresponding {\em truth approximation} $\Phi (\bmu^1)$ from Algorithm~\ref{alg1} to form a (one-dimensional) RB space $W_1 = \{\Phi (\bmu^1)\}$.  Next, we solve equation \eqref{eq:RBnlsys} to obtain a RB approximation $\widehat{\Phi}(\bmu)$ for each parameter in $\Xi_{\rm train}$ together with an error bound $\Delta_1(\bmu)$. The greedy choice for the $(i+1)$th parameter ($i = 1, \cdots, N-1$) is made by
\begin{equation}
 \bmu^{i+1} = \arg \max_{\bmu \in \Xi_{\rm train}}  \Delta_i (\bmu).
\end{equation}
The error bound is traditionally residual-based (i.e. $(F - \mathcal{L}_\mathcal{N}(\bmu) \widehat{\Phi}(\bmu))$-based) \cite{Rozza2008, Quarteroni2015, HesthavenRozzaStammBook}, 
\begin{equation}
\Delta_i (\bmu) = \frac{|| F -
\mathcal{L}_\mathcal{N}(\bmu) \widehat{\Phi}(\bmu)||_2} {\sqrt{\beta_{LB}(\bmu)}}, 
\end{equation}
where $\beta_{LB}(\bmu)$ is the smallest eigenvalue of $\mathcal{L}_\mathcal{N}(\bmu)^T \mathcal{L}_\mathcal{N}(\bmu)$. For simplicity, we compute it from the linear part of the operator \eqref{linearpb}. This stability constant is usually calculated by the Successive Constraint Method (SCM) \cite{HuynhSCM, CHMR-M2an, HKCHP} which is not necessary as the parameter dependence of the eigenvalue is obvious. 
After obtaining the $(i+1)$th {\em truth approximation} $\Phi(\bmu^{i+1})$ from Algorithm.~\ref{alg1}, we augment the  RB space by $W_{i+1}= \hbox{GS}(W_i, \vec{\Phi}(\bmu^{i+1}))$,
where $\mbox{GS}(W_{i},  \vec{\Phi}(\bmu^{i+1}))$ denotes a Gram-Schmidt orthogonalization procedure of the new vector $\vec{\Phi}(\bmu^{i+1})$ with respect to the previously selected and orthogonalized basis vectors in $W_{i}$. The detailed numerical scheme is described in Algorithm. \ref{algorithm0}.

\begin{algorithm}[h]
\caption{Reduced basis greedy sampling algorithm}
\label{algorithm0}
\vspace{0.5ex}
0. Choose  $\bmu^1$ randomly in $\Xi_{\rm train}$ and solve $\Phi (\bmu^1)$ from  Algorithm.~\ref{alg1}  \\[0.5ex]
1. Initialize $W_1 = \left\{\Phi (\bmu^1) \right\},$  \\[0.5ex]
2. \mbox{\textbf{For}} $N = 2,\ldots, N_{\max}$  \\[0.5ex]
3. $\quad\ \mbox{Solve }  c (\bmu) $ from Algorithm.~\ref{alg2}, and calculate $\Delta_{N-1}(\bmu)$  for all $\bmu \in \Xi_{\rm train}$  \\[0.5ex]
4. $\quad\ \mbox{Find } \bmu^{N} = \arg \max_{\bmu \in \Xi_{\rm train}}  \Delta_{N-1}(\bmu)$ \\[.5ex]
5. $\quad\ \mbox{Solve }  \Phi (\bmu^N)  ~ \mbox{from Algorithm}.~\ref{alg1}  ~ \mbox{and orthogonalize } W_{ N} = \mbox{GS}(W_{N-1},  \Phi (\bmu^N))$. \\[0.5ex]
6. \mbox{\textbf{End For}}
\end{algorithm}

With the greedy sampling algorithm, we still have to form and solve, in each iteration, the smaller RB systems in Step 3 for all $\bmu$ in the training set. While solving the RB system is inexpensive, forming it can be much more expensive. However, for a particular set of parametrized linear systems,  the RB system can be formed efficiently. The technique is an offline-online decomposition which is the topic of the next subsection.

\subsubsection{Offline-online computational procedure}
\label{sec:RBMcomputation}
To describe this procedure, we revisit the original FDM scheme \eqref{Operator} for the nonlinear PB equation. $\mathcal{L}_\N (\bmu)$ and $\vec{F}(\bmu)$ are the discretized operator and  the right hand side vector. For the purpose of forming the RB operator $A(\bmu)$ in equation \eqref{eq:RBnlsys} efficiently, we decompose $\mathcal{L}_\N (\bmu)$ at each $m$th iteration as follows {
\begin{equation}
\calL_\N(\bmu; \Phi_m) = D \calL^1_\N + \calL^2_\N(\bmu; \Phi_m) 
\end{equation}
Here, $D\calL^1_\N$ is the first part ($-D \nabla^2$) of equation \eqref{linearpb} having an explicit $\bmu$-dependence (note $\bmu = (D, V)$) but is $\Phi_m$-independent }
\begin{align}
\left(\mathcal{L}^1_\N \Phi_{m+1}\right)(x_j, y_k) &= - \frac{\Phi_{m+1}(x_{j-1}, y_k)-2\Phi_{m+1}(x_j, y_k) +\Phi_{m+1} (x_{j+1}, y_k)}{h_x^2} \nonumber\\
&- \frac{\Phi_{m+1} (x_j, y_{k-1}) -2\Phi_{m+1}(x_j, y_k) +\Phi_{m+1} (x_j, y_{k+1})}{h_y^2},  
\end{align}
while $\calL^2_\N(\bmu; \Phi_m)$ denotes the second part of equation \eqref{linearpb} that depends on both $\bmu$ and $\Phi_m$. When applied to $\Phi_{m + 1}$, it produces
$$\left(\calL^2_\N(\bmu; \Phi_m)\Phi_{m+1}\right)(x_j, y_k) = \cosh(\Phi_{m}(x_j, y_k)) \Phi_{m+1} (x_j, y_k).$$
Here $h_x =2/ N_x, h_y =2/ N_y$, with $N_x$ and $N_y$  being the the number of intervals in $x$-direction and $y$-direction, respectively.
After this decomposition, $A(\bmu; \widehat{\Phi}_m)$  can be written as
\begin{align}
A(\bmu; \widehat{\Phi}_m) = W_N^T \mathcal{L}_\N(\bmu; \widehat{\Phi}_m) W_N = A_1(\bmu) + A_2(\bmu; \widehat{\Phi}_m),
\end{align}
where
\begin{align*}
A_1(\bmu) & = D W_N^T \calL^1_\N W_N,\\
A_2(\bmu; \widehat{\Phi}_m) & = W_N^T \calL^2_\N(\bmu; \widehat{\Phi}_m) W_N.
\end{align*}

Below is a summary of the decomposition and the operation count each step takes.
\begin{itemize}
\item Realizing $W_N^T \calL^1_\N W_N$ is $\bmu$-independent, we precompute it by gradually populating this $N\times N$ matrix as we identify the RB bases in the space $W_N$ one by one. Indeed, when the $i$th basis is determined, we populate the $i$th row and $i$th column of $W_N^T \calL^1_\N W_N$. This step takes $O(N^2)$ operations.
\item Update $A_2(\bmu; \widehat{\Phi}_m) $ and the right hand vector
$W_N^T \vec{F}(\widehat{\Phi}_m)$ 
 at each iteration of the online procedure. This step takes $O(\N_0 N)$ operations.
\item Invert the RB matrix $A(\bmu)$ with $O(N^3)$ operations.
\item Form $\widehat{\Phi}(\bmu) = W_N c(\bmu)$ after each iteration, taking $O( \N_0 N)$ operations.
 \end{itemize}
Therefore, the total operation count of the online stage is
$$O( \N_0 N + N^3).$$
Although having $\N_0$-dependence which can be eliminated by the Empirical Interpolation Method \cite{Barrault2004, grepl2007efficient}, we note that the dependence is linear and it still produces an approximation much faster than the original FDM scheme.
From complexity analysis, this is possible because the coefficient matrix for the case with two-dimensional physical domain  in Algorithm.\eqref{alg1} is an $\N_0 \times \N_0$  banded matrix with a band width $2 N_x-1$.
The fact that $\N_0 = (N_x-1)(N_y+1)$ is often very large and  $N$ is  typically very small and in particular $N < N_x$ are an indicator that our RB algorithm will be much faster than the truth solver. This is indeed corroborated in the next section where we present numerical results.

\section{Numerical examples}
\label{sec:result}

In this section, we show numerical results in both one- and two-dimensional spaces to demonstrate the performance of the proposed RBM algorithm. The common parameter space $\D \ni (D, V)$ is taken to be $\D = [0.08^2,0.4^2] \times [0, 5]$,
which is discretized by a so-called {\em training set}
\[
\Xi = \{D | \sqrt{D}  \in (0.08:0.02:0.4)\} \times \{(0:0.25:5)\}. 
\]
We also define a testing set
\[
\Xi_{\rm test} = \{D | \sqrt{D} \in (0.085:0.01:0.395)\} \times \{(0.4:0.5:4.4)\} 
\]
which in particular does not {intersect} with the training set. 
Here, the notation $a : h : b$ denotes an equidistant discretization of the interval $[a, b]$ by elements of size $h$. 
In addition, we let $E(N)$ represent the maximum relative error over all $\bmu$ in $\Xi_{\rm test}$ of the reduced basis solution using $N$ basis functions in comparison to the truth approximation
\begin{equation}
E(N) =  \max_{\bmu \in \Xi_{\rm test}}\{\| \Phi(\bmu) - \widehat{\Phi}(\bmu)\|_\infty/\|\Phi\|_{L^\infty(\Xi_{\rm test}, L^\infty(\Omega))}\}
\end{equation}
where 
\[
\|\Phi\|_{L^\infty(\Xi_{\rm test}, L^\infty(\Omega))} = \max_{\bmu \in \Xi_{\rm test}} \|\Phi(\bmu) \|_\infty.
\]
Lastly, we let 
\[
\Delta_{\rm RB}^{\max}(N) =\max_{\bmu \in \Xi} \Delta_N(\bmu)
\]
represent  the maximum error bound over the discretized set $\Xi$ when $N$ parameters are selected. 

\subsection{One dimensional  space}
\label{sec:onedim}
We consider the PB equation \eqref{eq:PB} with source $g(x,y)=0$. The problem reduces to one-dimensional thanks to the homogeneity in the $y$ direction. The truth approximation $ \Phi$ is obtained with a central finite difference scheme in Algorithm.~\ref{alg1}.
The physical domain $[-1, 1]$ is divided into
$N_x$ cells by $\mathcal{N} = N_x+1$ grid points. 

\begin{figure}[htp!]
\includegraphics[scale=0.32]{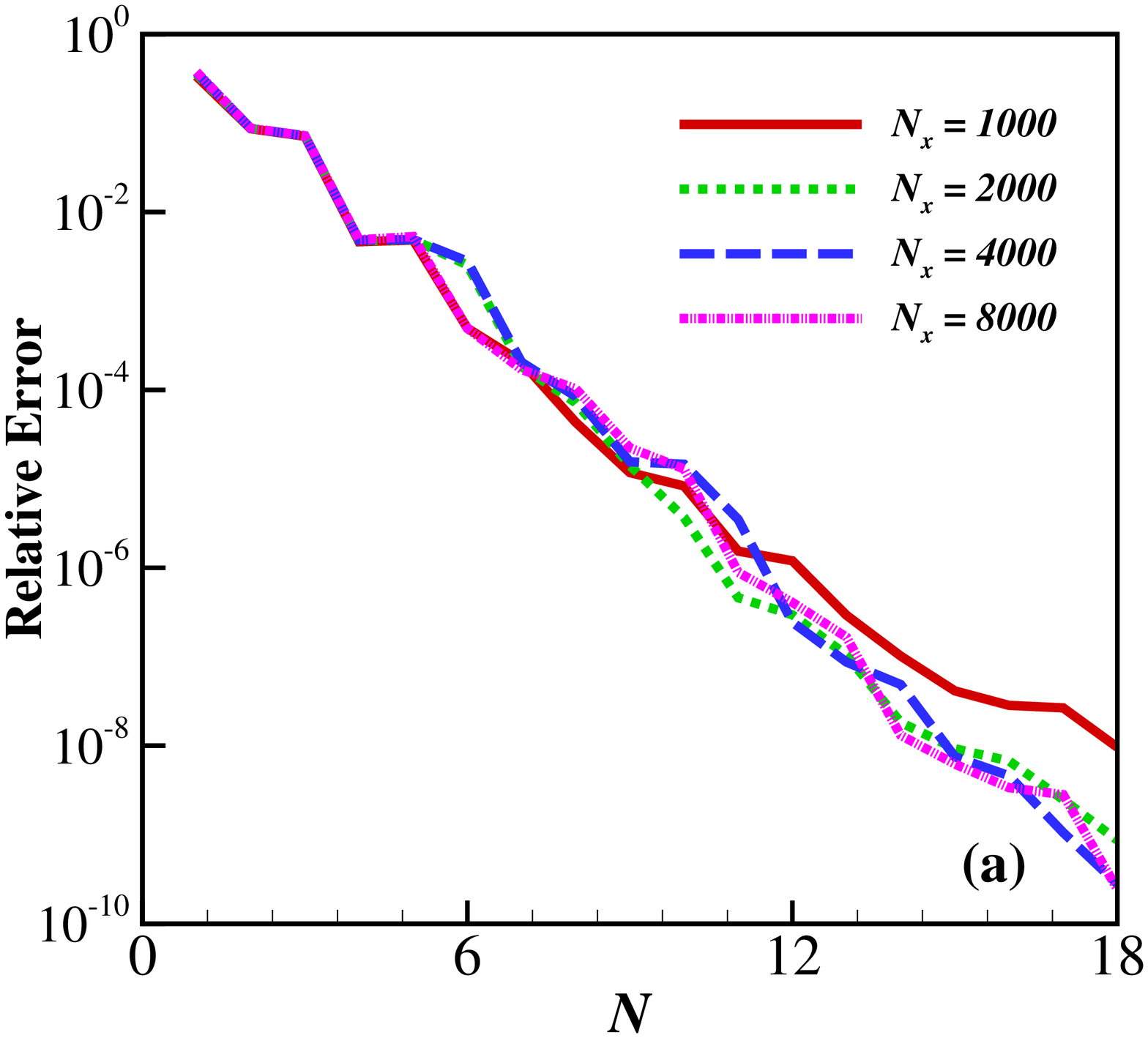}
\includegraphics[scale=0.32]{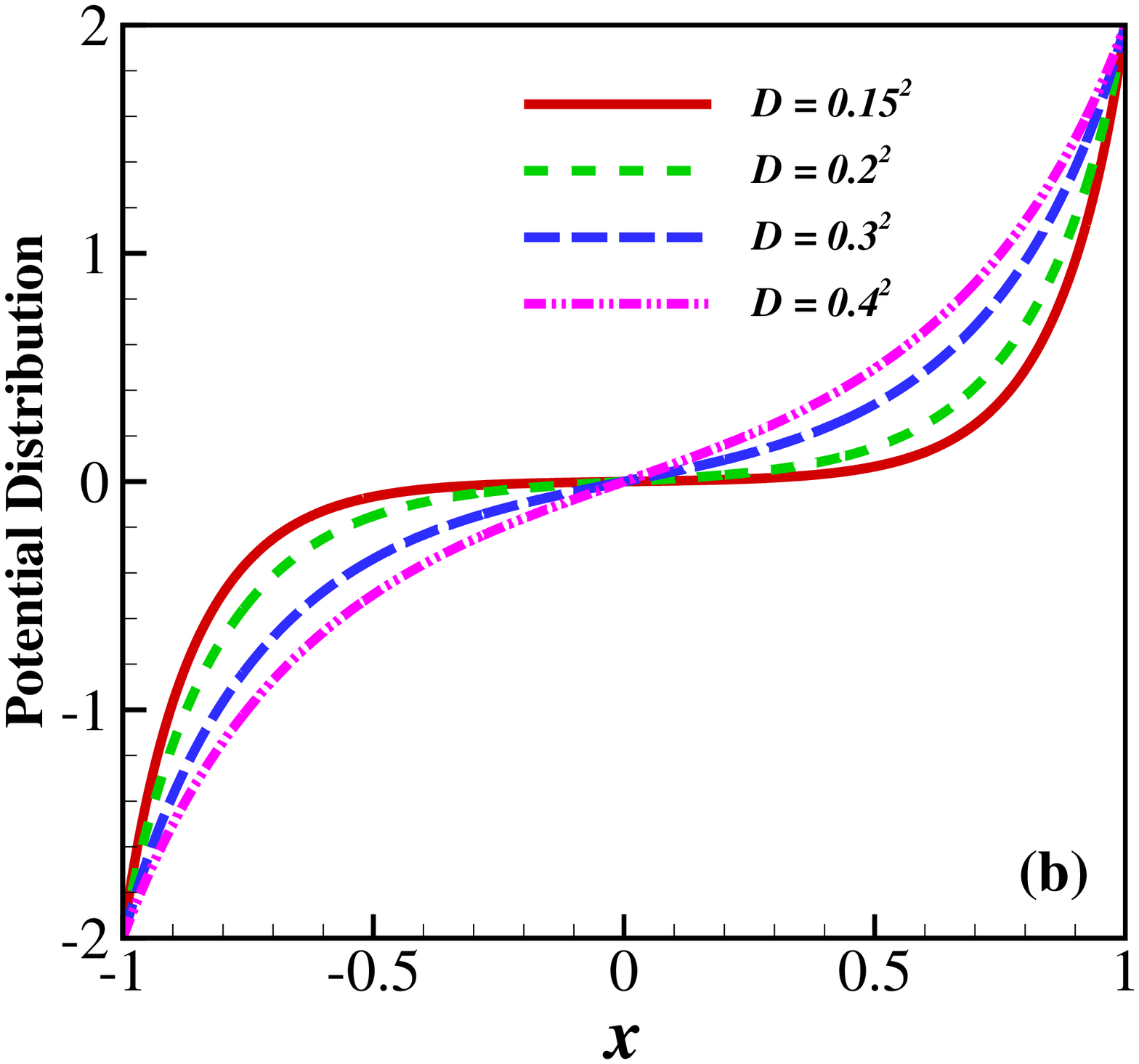}\\
\includegraphics[scale=0.32]{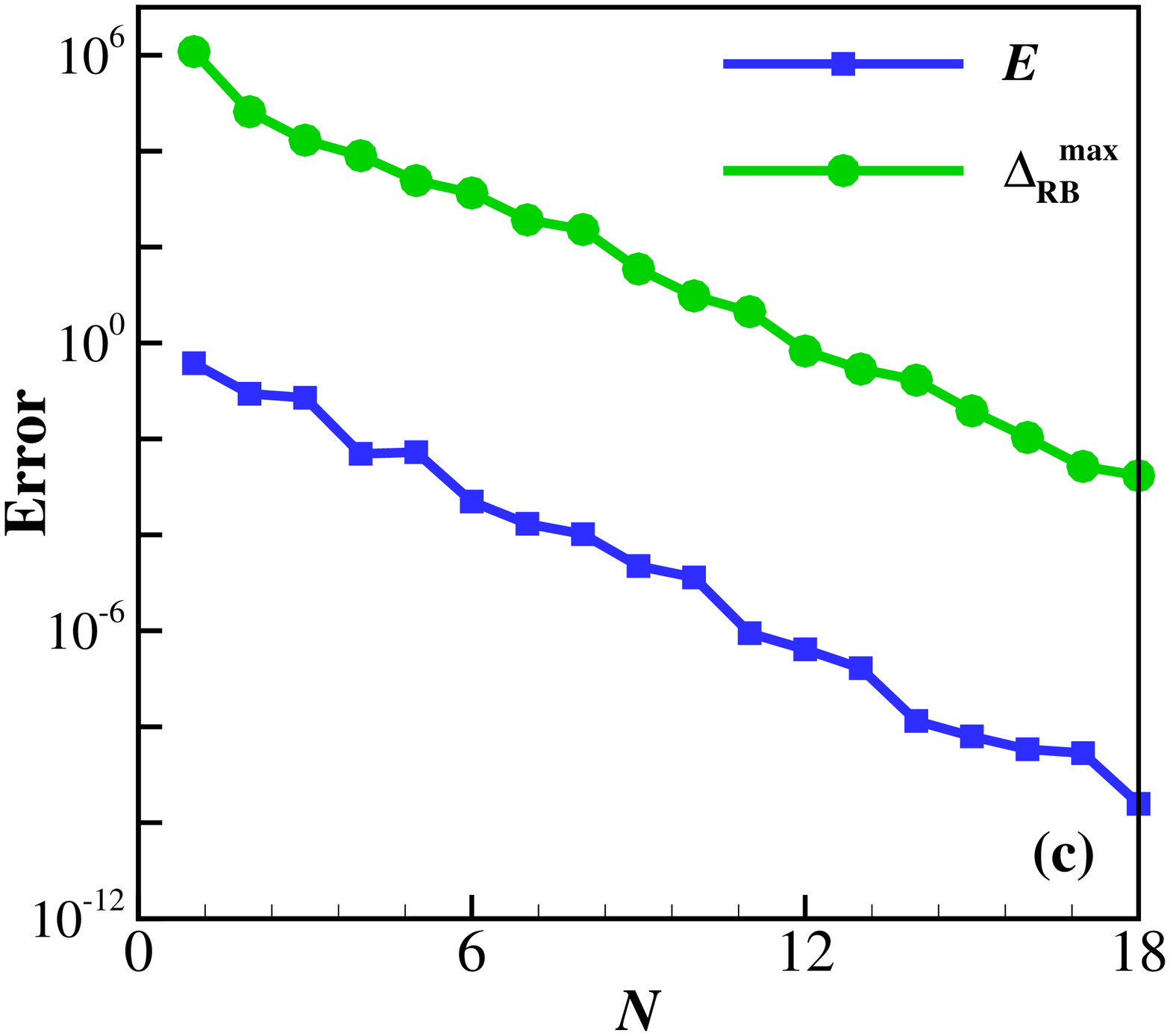}
\includegraphics[scale=0.32]{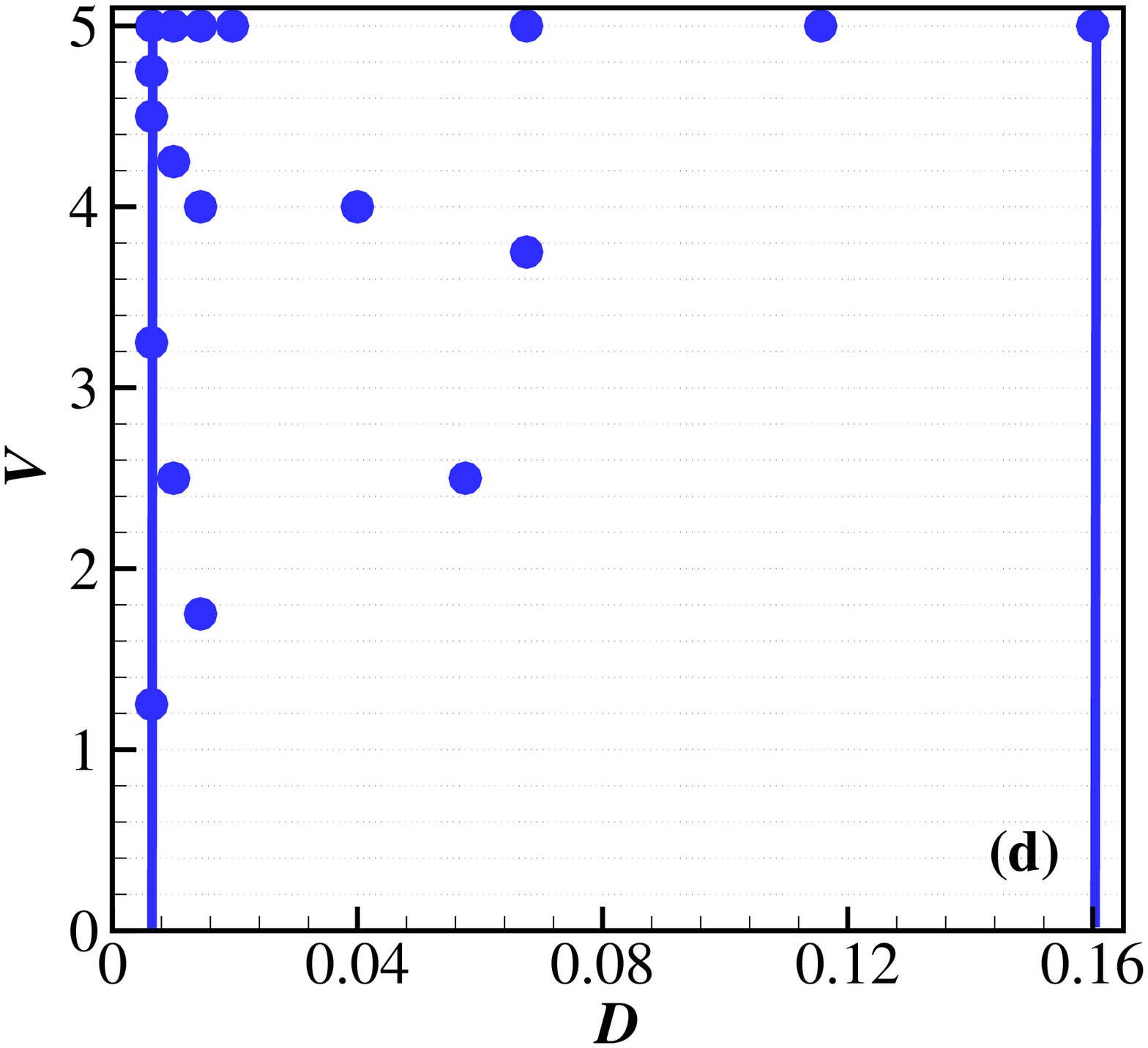}
\caption{(a) The maximum relative errors of RB approximations $E(N)$ versus $N$. (b) Potential distributions at representative $D$ values with $V = 2$. (c) Maximum error estimators  of RB approximation $\Delta_{\rm RB}^{\max}$  versus $N$.  (d) The location of the selected parameters where the reduced bases are computed.}
\label{FIG1}
\end{figure}

Fig. \ref{FIG1}(a) displays the relative errors of the RB solution when different partition numbers $N_x =1000, 2000, 4000$ and $8000$ are used for the truth approximation.
One clearly observes exponential convergence of the error with respect to the number of reduced bases.
For this one-dimensional problem, using $N=12$ basis functions is enough for the error $E(N)$ to reach $\sim10^{-6}$ which is of the same magnitude as the truth approximation error. Fig. \ref{FIG1}(b) shows the sample solutions $\Phi(x; \bmu=(D, 2))$, i.e. the potential distributions at different $D$ values with $V =2$ and $N_x =10,000$. 
It is clear that the smaller $D$ is, the stronger the boundary-layer is, a manifestation of the
nonlinearity of the PB equation. 
This provides an intuitive account of why the parameter locations for the chosen RB snapshots, shown in Fig.   \ref{FIG1}(d), clusters around smaller $D$ values.

Now we fix $N_x = 10,000$ and study the effectivity of the RB model.
We show the comparison of $\Delta^{\max}_{RB}(N)$ and  the RB relative error $E(N)$ in Fig.~\ref{FIG1}(c),
and the selected parameters' locations  in Fig.~\ref{FIG1}(d). 
It is noted that the error estimator is decreasing with similar exponential speed as the error even though the error is calculated in a stronger norm ($L^\infty$) than that for the residual in the error estimator. The distribution of the chosen parameters shown
in Fig.~\ref{FIG1}(d) accurately reflects the property of the nonlinear PB equation, i.e.,  most parameters with smaller $D$ value are selected and less are chosen with large $D$ value. This is because that the elliptic PDE with a small $D$ tends to form a boundary layer near the boundary.
In addition, most of the selected parameters are on the boundaries, i.e. when $V = 0$ or  $5$ and $\sqrt{D} = 0.08$ or $0.4$, the latter denoted by the left ($D = 0.0064$) and the right ($D = 0.16$) vertical lines.

\begin{figure}[htp!]
\includegraphics[scale=0.32]{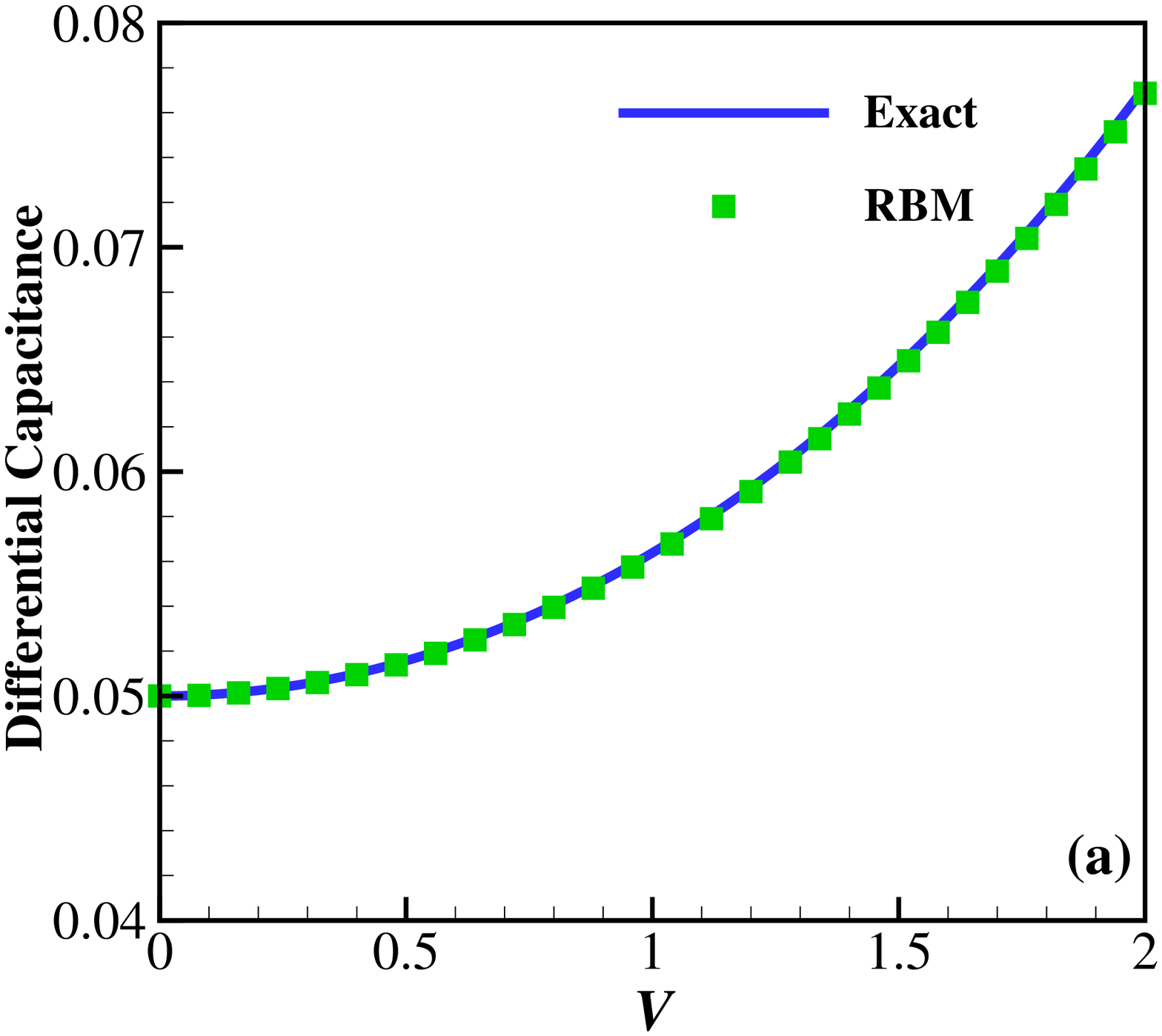}
\includegraphics[scale=0.32]{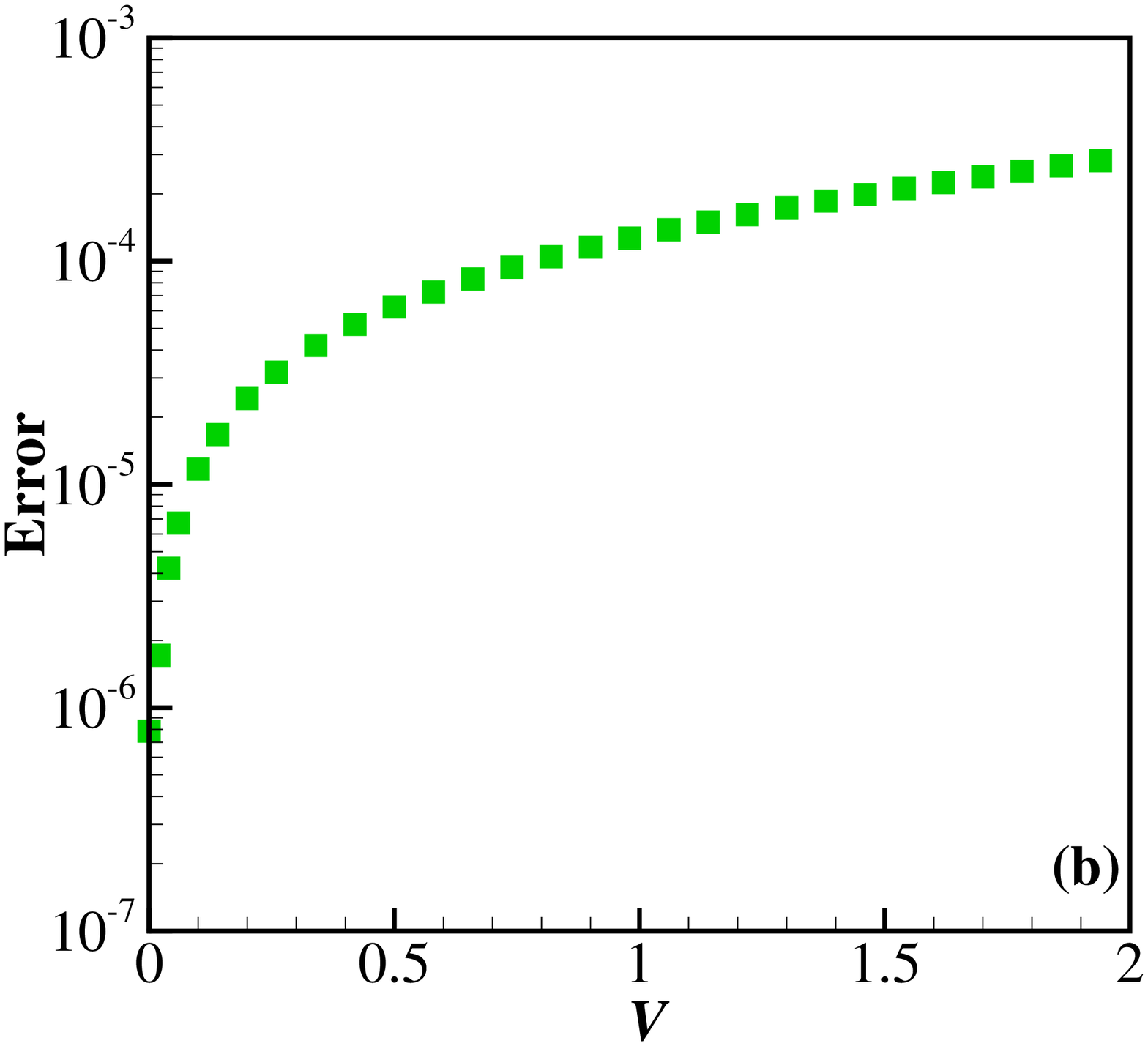}
\caption{(a) Differential capacitance $C$ as function of the voltage $V$. (b)  The error of the differential capacitance as function of $V$.}
\label{FIG2}
\end{figure}

Next, we report the result on calculating the total differential capacitance for the symmetric electrolyte. 
In the FDM scheme, the numerical surface charge density $\sigma$ is calculated by,
\[
\sigma =  D \partial_x \Phi(x = -1)   = D \frac{4 \Phi(x_2)- 3 \Phi(x_1)-\Phi(x_3)} {2h}
\]
upon which the differential capacitance is obtained.
We take $D = 0.01$, $ V = 0:0.02:2$ and $N =16$ to calculate a surrogate differential capacitance.  Fig. \ref{FIG2}(a)
shows the results of the  the capacitance as function of the boundary
voltage $V$ by the RBM, compared to the exact curve \eqref{eq:capacitance}, and panel (b) shows the error.
It can be observed that the RBM is very accurate even with a small number of bases ($N=16$).

\subsection{ Two dimensional space}
\label{sec:twodem}

In this section, we solve a two dimensional PB equation \eqref{eq:PB} with the {fixed} charge $g(x,y)= \exp[-50(x^2+y^2)]$.
The  total number of nodes for the central-difference scheme in Algorithm.~\eqref{alg1}  is $\mathcal{N} = (N_x+1)  (N_y+1)$,
and the number of  unknown nodes is $\mathcal{N}_0 = (N_x-1)  (N_y+1)$.
The boundary conditions on $y =\pm 1$ are approximated by the central difference scheme.

\begin{figure}[htp!]
\includegraphics[scale=0.45]{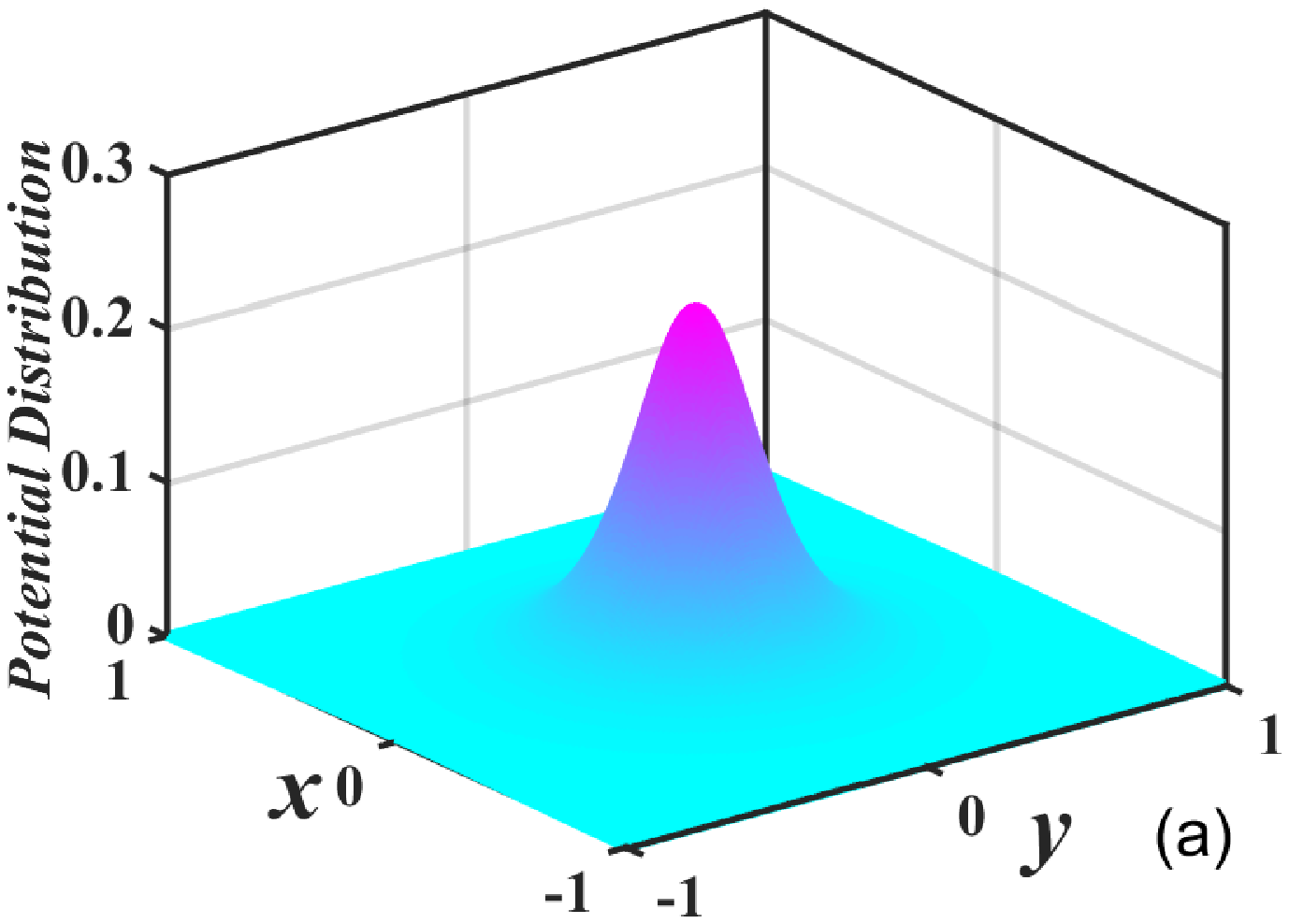}
\includegraphics[scale=0.45]{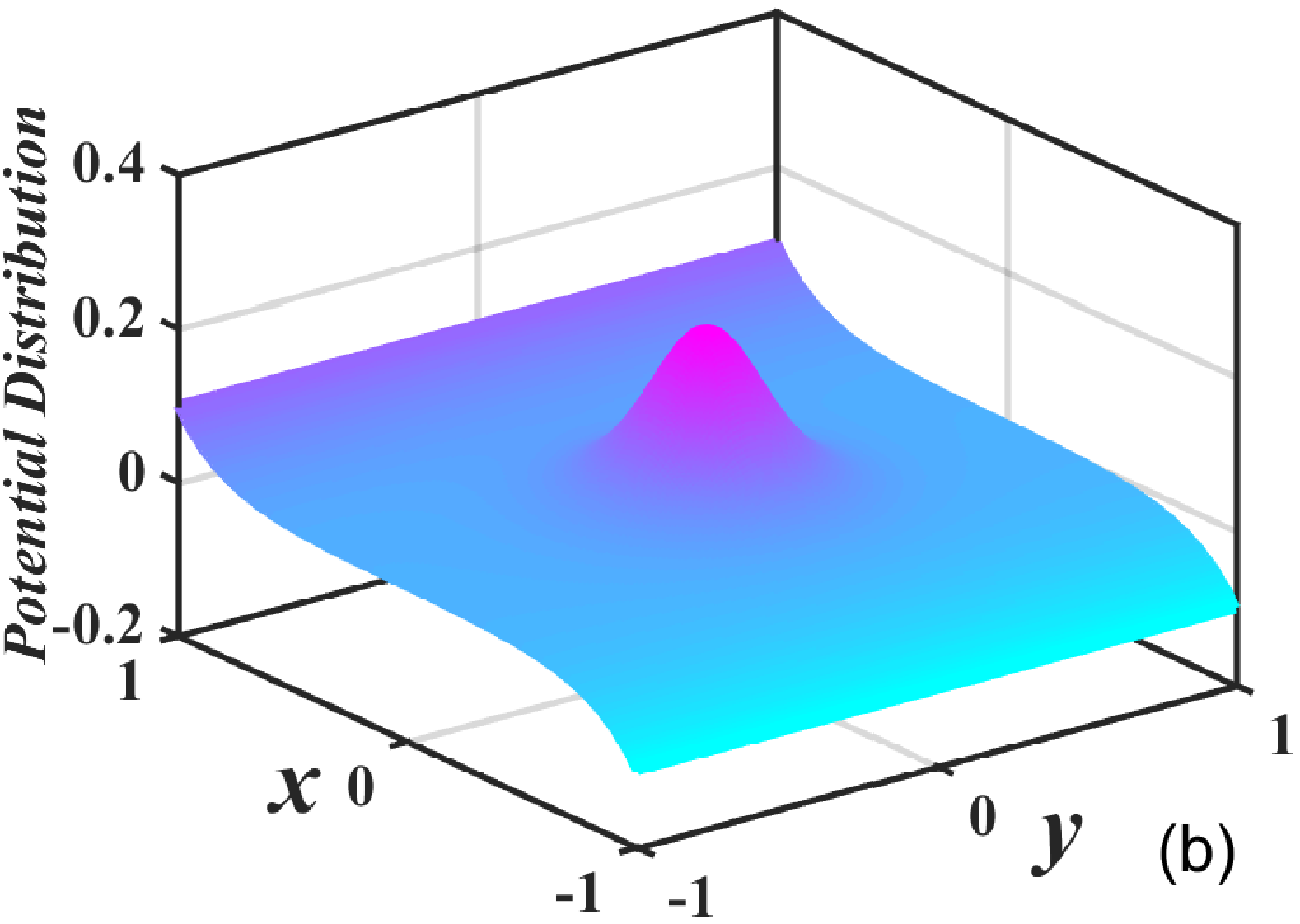} \\
\includegraphics[scale=0.45]{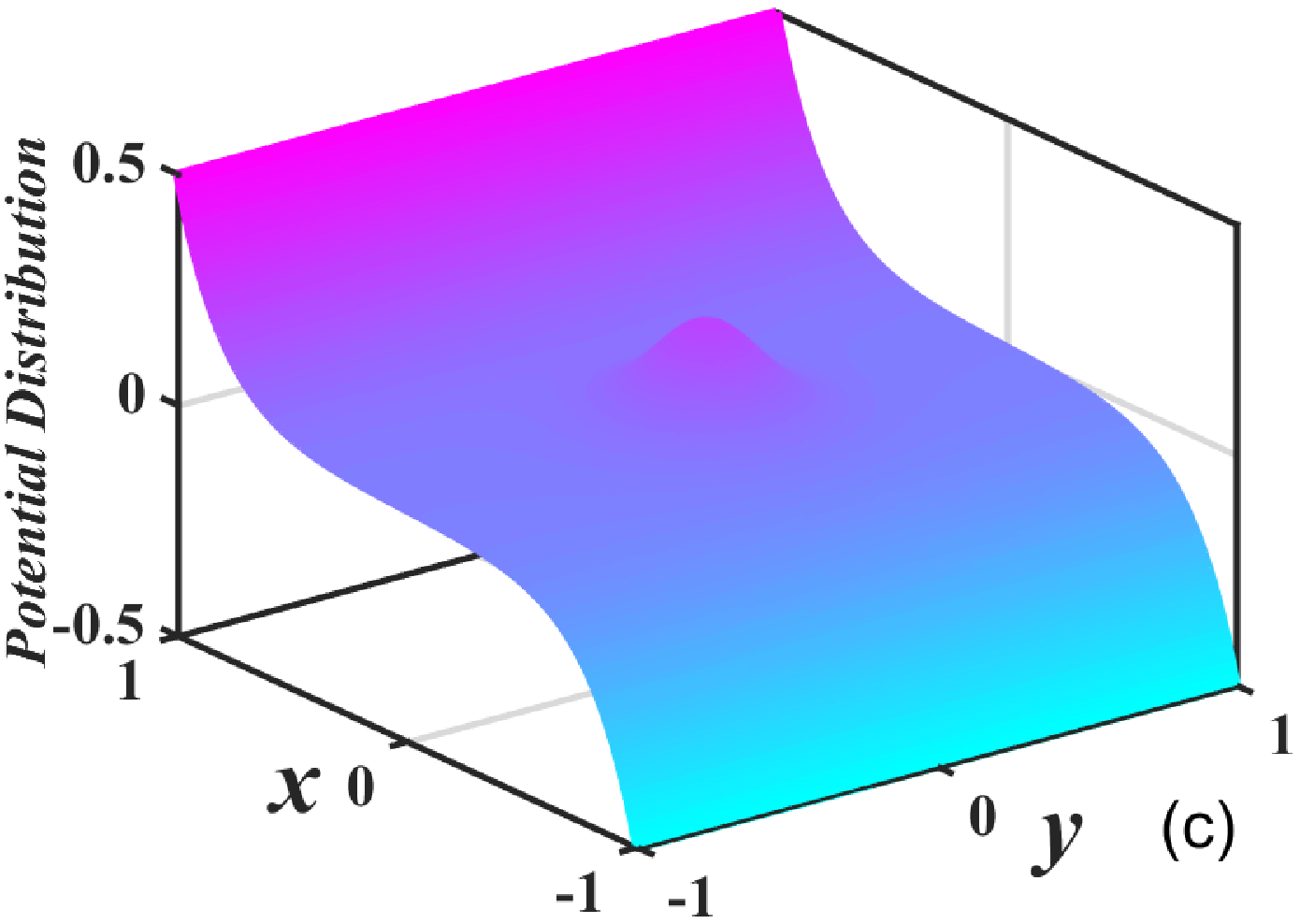}
\includegraphics[scale=0.45]{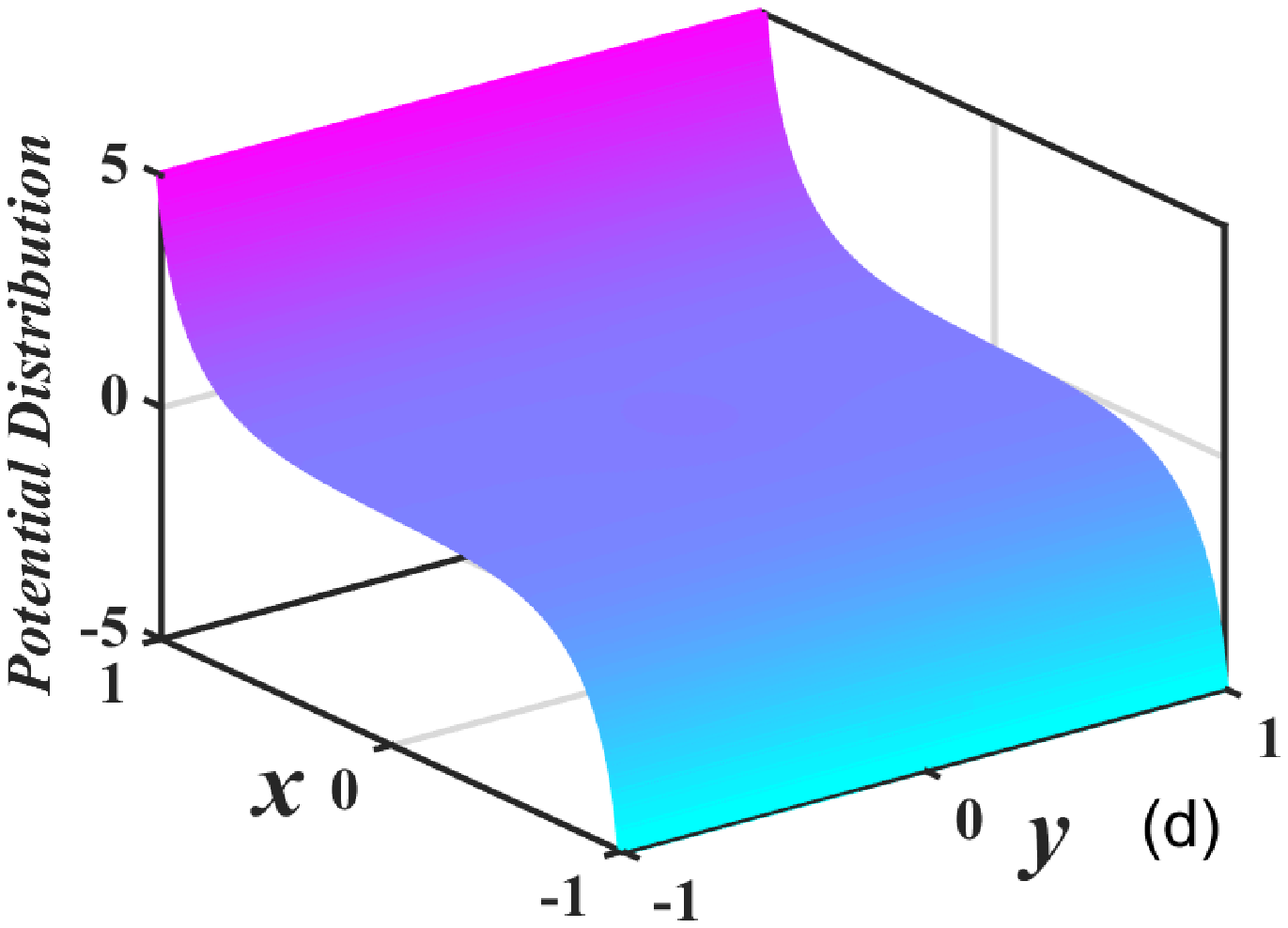}
\caption{Potential distributions for different $V$ by the RBM. (a) $V = 0$ (b) $V= 0.1$ (c) $V= 0.5$ (d) $V = 5$.}
\label{FIG3}
\end{figure}

\begin{figure}[htp!]
\includegraphics[scale=0.32]{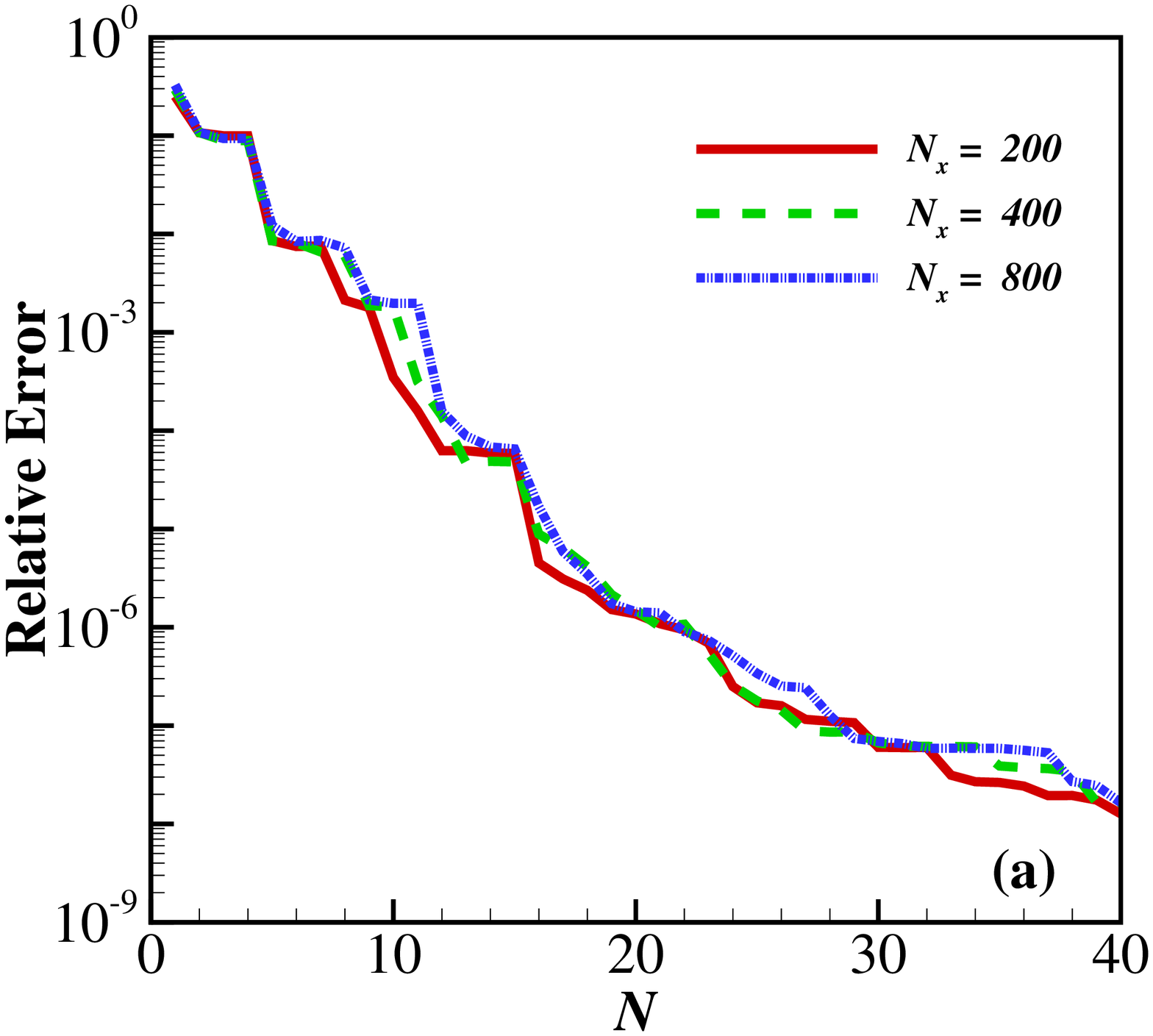}
\includegraphics[scale=0.32]{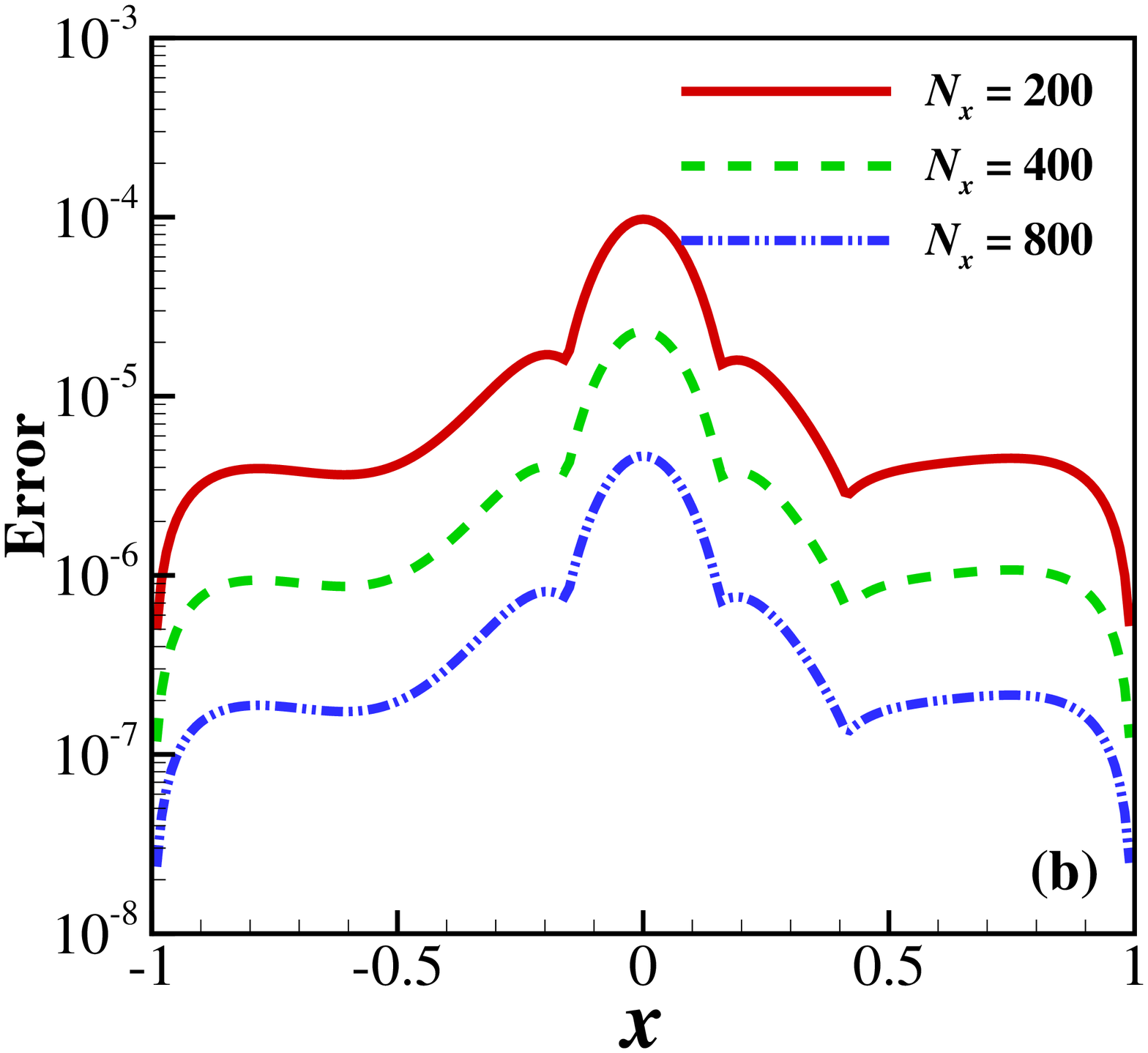} \\
\includegraphics[scale=0.32]{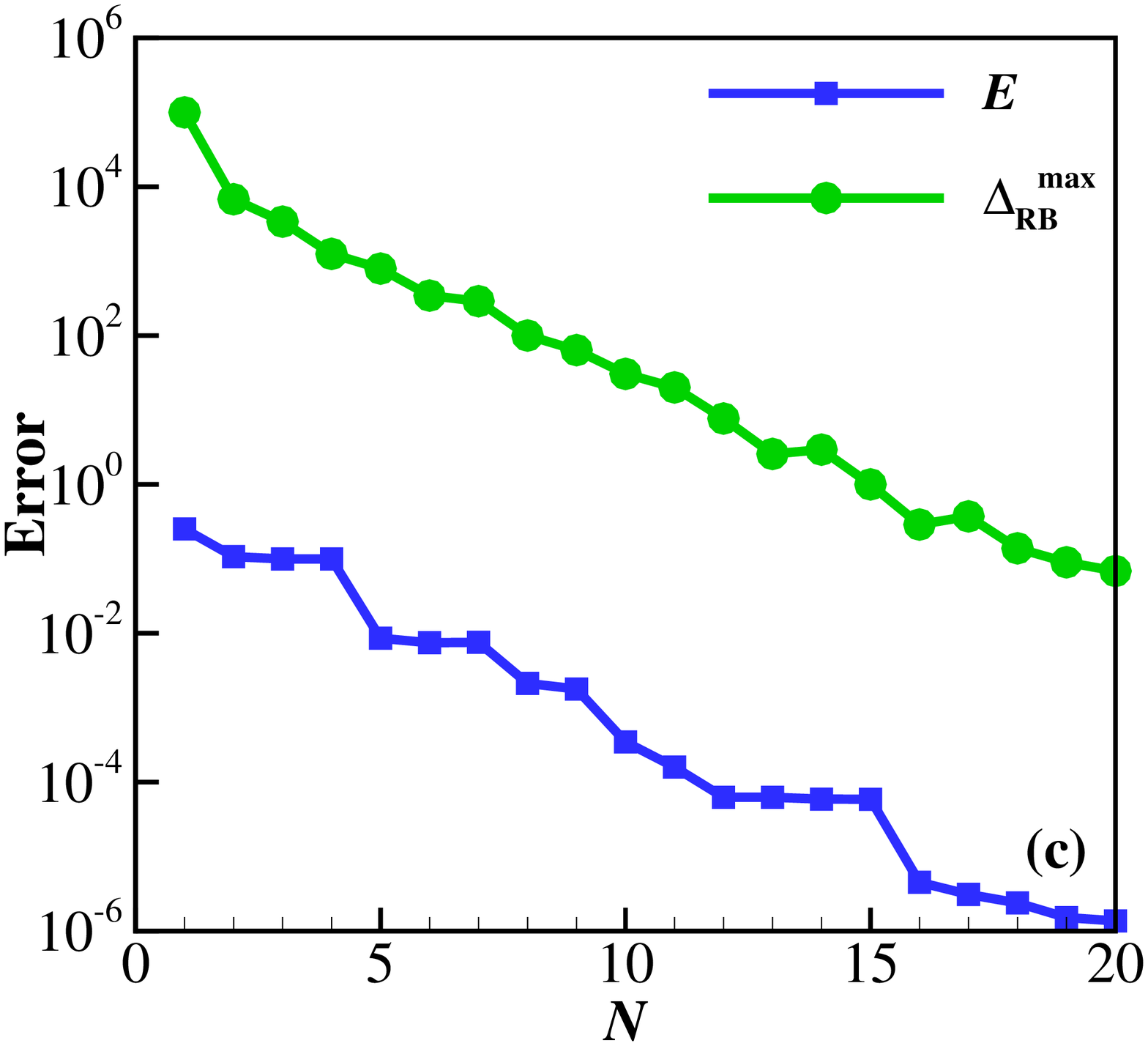}
\includegraphics[scale=0.32]{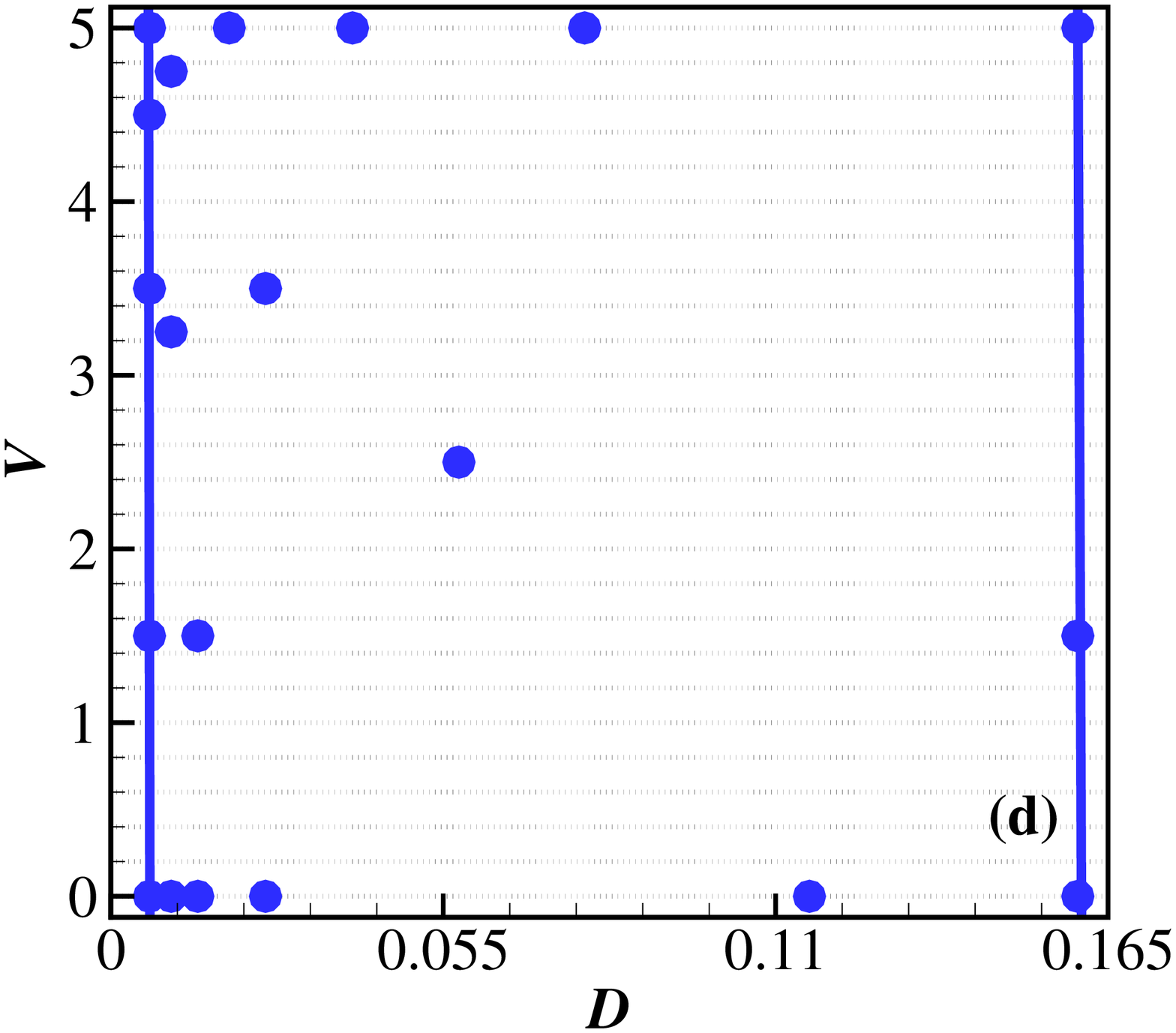}
\caption{(a) Convergence of RB approximation of different partition numbers at different basis numbers $N$. \rev{(b) The error accuracy $	E_x$ of FDM in x-direction at different partition numbers $N_x$ with $V = 0.1, D = 0.04$.} (c) Maximum error estimators $\Delta_{\rm RB}^{\max}$  at different basis numbers $N$. (d) Selected parameters' location.}
\label{FIG4}
\end{figure}

We first show representative RB approximations $\widehat{\Phi}$  when $D = 0.04, N_x = N_y =200$ in Fig.~\ref{FIG3}.   It can be observed that when the applied voltage is smaller, the variation of the potential distribution $\Phi$ is stronger.
Using the same parameter subsets as those in the previous section, we show the RB relative error $E(N)$ in Fig.~\ref{FIG4}(a) and we again observe the convergence for all partition numbers $N_x, N_y$ with $N_x = N_y =200,  400, 800$. We see that, for the RBM solution to approximate the FDM sufficiently closely, a basis number $N$ larger than the one-dimensional case is necessary. In order to show that when $N \approx 20$ we are at the FDM accuracy level, we verify our FDM accuracy in Fig.~\ref{FIG4}(b). Indeed, we set $N_x = 200, 400, 800$ and $1,600$ with $V = 0.1, D = 0.04$ and take the solutions with $N_x = 1,600$ to be the reference, and define $$E_x(i) = ||\Phi_{N_x}(i,\cdot) - \Phi_{1600}(i,\cdot)||_\infty , i = 1,2,\cdots,N_x+1.$$
 $E_x$ can be viewed as the absolute error at each discrete  node in x-direction and the infinite norm is for the y-direction.
The distribution of $E_x$ is  shown in Fig.~\ref{FIG4}(b) and clearly we only need about 20 parameters to reach the truth relative error $\sim 10^{-6}$ which is the accuracy level of the FDM solution. 
The comparison between $E(N)$ and $ \Delta_{RB}^{\max}$ when { $N_x = N_y =200, N =20$ is shown in Fig.~\ref{FIG4}(c). Parameter locations for the RB snapshots are shown  in Fig.~\ref{FIG4}(d) when $N_x = N_y =200, N =20$} and we also find that small voltages are often selected,  which is consistent with the conclusions in
Fig.~\ref{FIG3}. Specifically,  most of selected parameters sit on the boundaries of the parameter domain, and the left vertical line in Fig.~\ref{FIG4}(d) refers to $D = 0.0064$, and the right one refers to $D = 0.16$.
Time consumption comparison is shown in Table.~\ref{time2}. Obviously the online time of  RB approximation is smaller than Algorithm.~\eqref{alg1}. This is mainly due to the coefficient matrix structure of the two dimensional case. It is not a tridiagonal matrix as in the one dimensional case and order reduction is significant for such matrix with a wide band width.

\begin{table}[htpb!]
	\centering
		\begin{tabular}{lrrrr}
		\hline
$N_x, N_y$& RBM  $\widehat{\Phi}$&  FDM $\Phi$   \\ \hline
100 	&  0.0231827  &    0.2228043  \\ \hline
200	&  0.0960147  &    1.1668299 \\ \hline
400	&   0.4744200     &  5.3805034  \\ \hline
800	&  2.0143542     &    27.2314357 \\ \hline
	\end{tabular}
	\caption{Online computational times  at different partition numbers  $N_x, N_y $  when $V = 3.57, D = 0.25^2, N = 20$.}
	\label{time2}
\end{table}

\begin{figure}[htp!]
\includegraphics[scale=0.32]{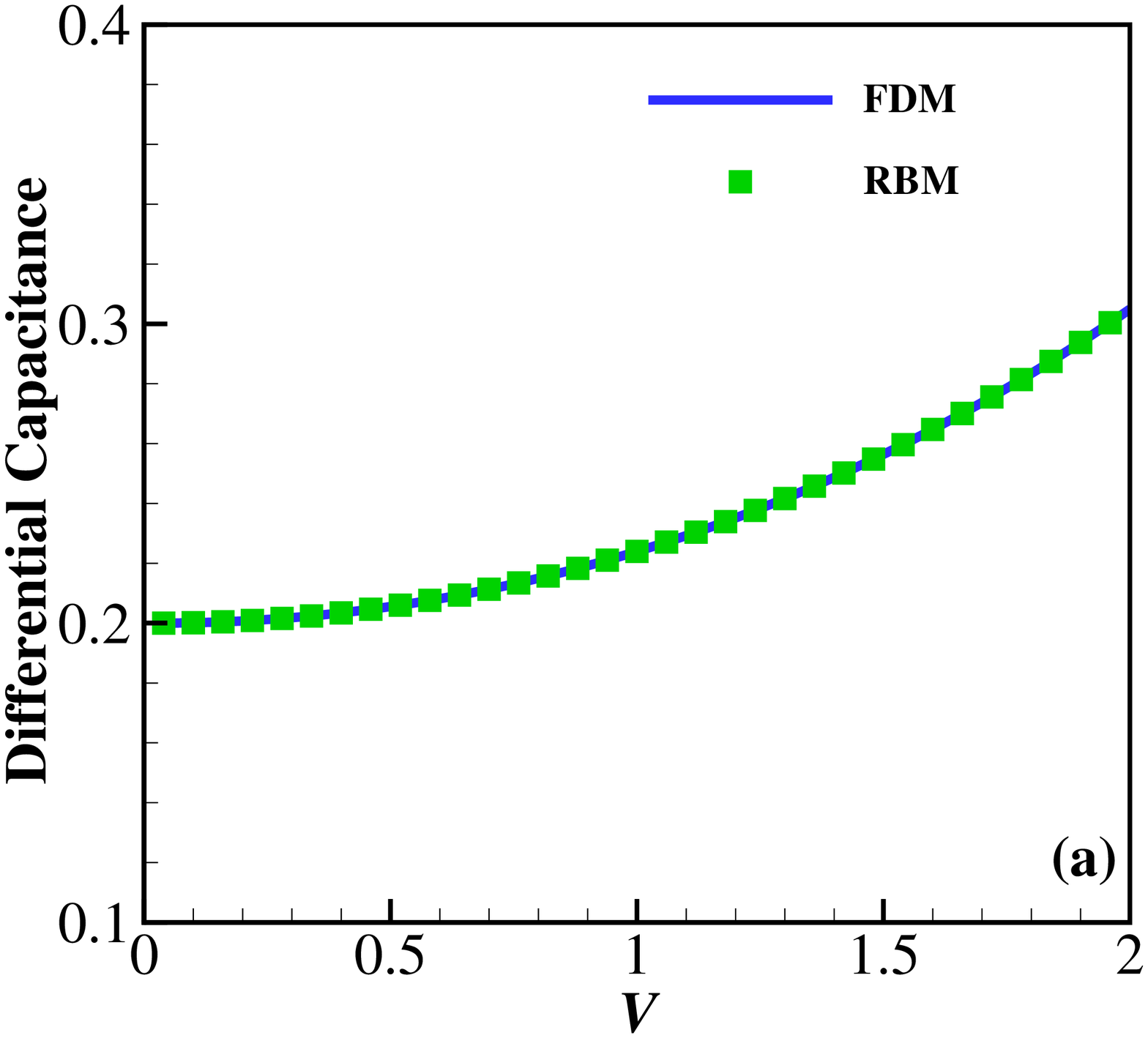}
\includegraphics[scale=0.32]{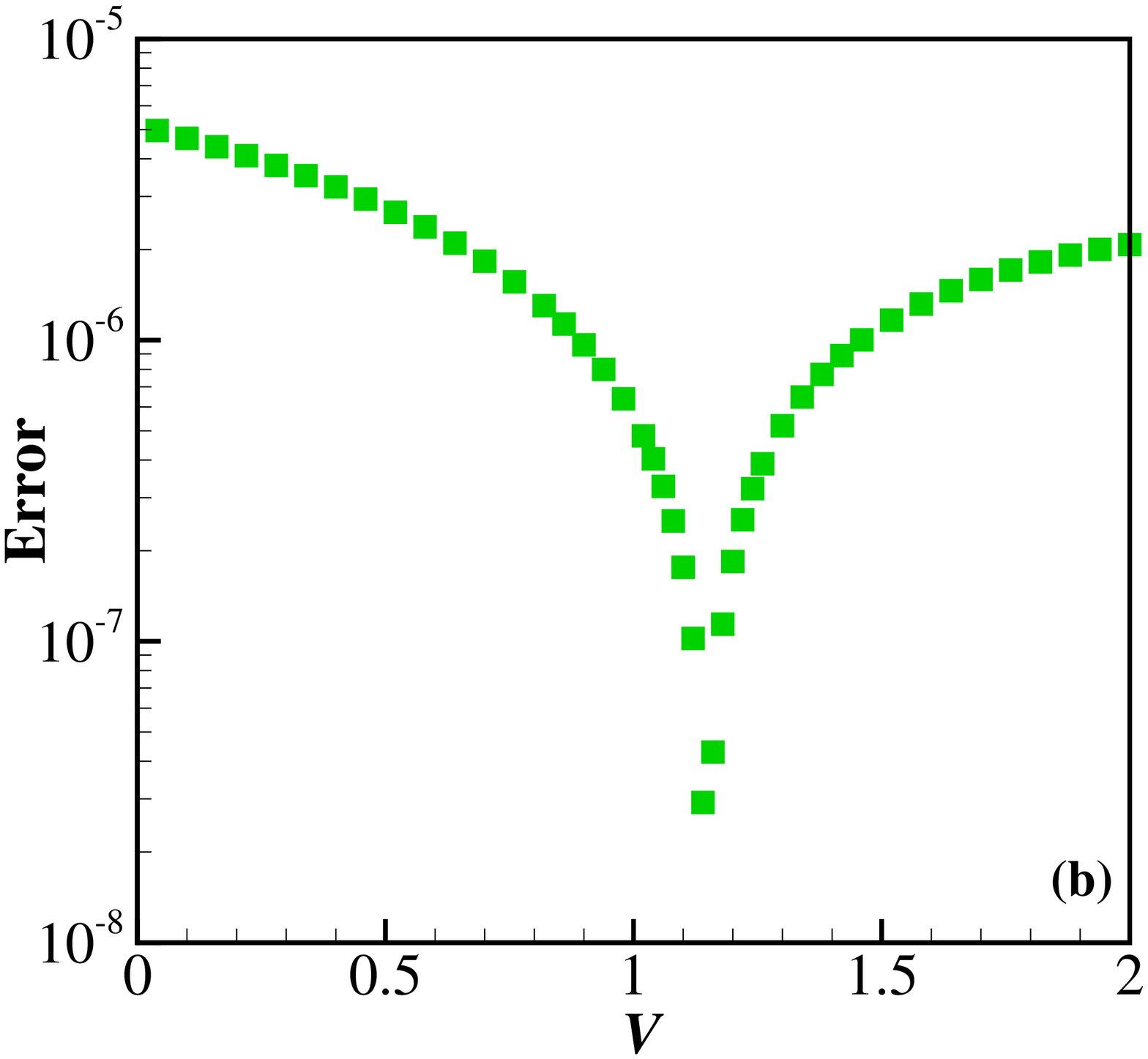}
\caption{(a) Differential capacitance $C$ at different voltages $V$. (b)  Corresponding error  at different voltages $V$.}
\label{FIG5}
\end{figure}

We take $N_x=N_y=800, D =0.04 $, $V = 0.04:0.02:2$ and $N= 20$ for comparison between the RB capacitance and differential capacitance derived by FDM. These results are showed in Fig.~\ref{FIG5}, calculating the differential capacitance $C$ for the two dimensional problem and its error. 
In fact, calculating these differential capacitance through the RBM is $15$ ore more times faster than using FDM.
This demonstrates that the RBM approximation is efficient in calculating the parameterized physical quantities such as the differential capacitance of the electrochemical system.

\section{Conclusion}
\label{sec:conclusion}
This paper applies the  RB algorithm to solve the parametrized nonlinear PB equation in both one and two dimensional physical spaces with a two dimensional parameter space. Though PB equation is non-affine and has exponential nonlinearity, our algorithm shows good accuracy and the selected parameters' distribution accurately reflects the nonlinearity of the PB equation. In future work, we consider further enhancement to the algorithm including application of EIM achieving total $\N$-independence of the online solver, and a novel approach that achieves the same efficiency without EIM.

\section*{Acknowledgments}
L. Ji and Z. Xu acknowledge the support from grants NSFC 11571236 and 21773165 and HPC center of Shanghai Jiao Tong University. Y. Chen was partially supported by the United States National Science Foundation grant DMS-1719698.



\end{document}